\documentclass[11pt,twoside,a4paper]{amsart}

\title{Essential representations of $C^*$-correspondences}
\author{Ilan Hirshberg}
\address{Department of Mathematics and Computer Science, University of Southern Denmark, Campusvej 55, DK-5230 Odense M, Denmark}
\email{ilan@imada.sdu.dk}

\usepackage{amssymb}
\usepackage{amscd}
\usepackage{amsthm}
\usepackage{xypic}

\theoremstyle{plain}

\newtheorem{Thm}{Theorem}[section]
\newtheorem{Cor}[Thm]{Corollary}
\newtheorem{Lemma}[Thm]{Lemma}
\newtheorem{Claim}[Thm]{Claim}

\newtheorem{Prop}[Thm]{Proposition}

\theoremstyle{definition}
\newtheorem{Def}[Thm]{Definition}

\newtheorem{Exl}[Thm]{Example}
\newtheorem{Rmk}[Thm]{Remark}
\newtheorem{Question}[Thm]{Question}

\newcommand{\B}{\mathcal{B}}
\newcommand{\A}{\mathcal{A}}
\newcommand{\M}{\mathcal{M}}
\newcommand{\J}{\mathcal{J}}
\newcommand{\K}{\mathcal{K}}

\newcommand{\T}{{\mathbb T}}

\newcommand{\Z}{{\mathbb Z}}
\newcommand{\C}{{\mathbb C}}
\newcommand{\EE}{\mathcal{E}}

\newcommand{\lb}{\left <}
\newcommand{\rb}{\right >}

\renewcommand{\arrow}{\rightarrow}

\newcommand{\TE}{\mathcal{T}_E}
\newcommand{\HTE}{\widehat{\TE}}

\newcommand{\OO}{\mathcal{O}}

\newcommand{\Jc}{\J_c}

\begin{document}

\begin{abstract}
Let $E$ be a $C^*$-correspondence over a $C^*$-algebra $\A$ with non-degenerate faithful left action. We show that $E$ admits sufficiently many essential representations (i.e. representations $\psi$ such that $\overline{\psi(E)H} = H$) to recover the Cuntz-Pimsner algebra $\OO_E$.   
\end{abstract}
\maketitle
\section{Introduction}

Let $\A$ be a $C^*$-algebra. Recall that by a $C^*$-correspondence over $\A$ we mean a Hilbert module $E$ over $\A$, along with a fixed left action of $\A$ on $E$ via adjointable operators (i.e. a homomorphism $\A \arrow \B(E)$ -- we suppress notation for this homomorphism). $\TE$ and $\OO_E$ denote the Pimsner-Toeplitz algebra and the Cuntz-Pimsner algebra, respectively (see \cite{pimsner}). We shall frequently identify $\A$ and $E$ with their images in $\TE$ and $\OO_E$.
 
\emph{Throughtout this paper, we shall assume that $E$ is full (i.e. $\overline{\lb E,E \rb} = \A$\footnote{We follow the convention that when we place two spaces into a bilinear or sesquilinear expression such as $\lb E,F\rb$, or $\A \B$, what is meant is the span of all possible pairings.}
), and that the left action of $\A$ on $E$ is faithful and non-degenerate}. 

By a \emph{representation} of $E$ on $H$ (c.f. \cite{ms1,FMR,me JFA}) we mean a pair $\psi_E : E \arrow \B$ (a linear map), $\psi_{\A}: \A \arrow \B$ (a homomorphism)
such that for all $e,f\in E$, $a\in \A$ we have $\psi_E(ea) = 
\psi_E(e)\psi_{\A}(a)$, $\psi_E(ae) = 
\psi_{\A}(a)\psi_E(e)$, $\psi_E(e)^*\psi_E(f) = \psi_{\A}(\lb 
e,f\rb)$. We will write $\psi$ for both $\psi_{\A}$ and $\psi_E$, when 
it causes no confusion. Note that what we call `representation' is called `isometric covariant representation' in \cite{ms1}. A representation $\psi$ of $E$ induces a representation of $E^{\otimes n}$, which we denote by $\psi$ as well, by $\psi(e_1\otimes \cdots \otimes e_n) = \psi(e_1)\cdots\psi(e_n)$.

Let $\EE = \bigoplus_{n=0}^{\infty}E^{\otimes n}$ (the Fock module). Each representation $\pi: \A \arrow \B(H)$ induces a representation of $E$ on $\EE \otimes_{\pi} H$ (where we use $\pi$ to view $H$ as an $\A-\C$ correspondence). Following Arveson's terminology for product systems (c.f. \cite{arveson's book}), we call a representation \emph{singular} if $\bigcap_{n=0}^{\infty} \overline{\psi(E^{\otimes n})H} = 0$, and \emph{essential} if $\overline{\psi(E)H} = H$ (essential representations are called `fully coisometric' in \cite{ms1}). 

Recall the following generalization of the Wold decomposition from \cite{ms1}.
\begin{Thm} Let $\psi$ be a representation of $E$ on $H$. $H$ decomposes uniquely to a direct sum of invariant subspaces $H_s \oplus H_e$, such that $\psi|_{H_s}$ is singular and $\psi|_{H_e}$ is essential. Furthermore, the singular part of the representation is unitarily equivalent to a representation induced from the Fock module (as explained above).
\end{Thm}

Thus is it easy to construct singular representations, however there doesn't seem to be an obvious way in general to construct essential ones. Any representation of $E$ gives rise to a representation of the Pimsner-Toeplitz algebra $\TE$ (see \cite{pimsner}, or \cite{me JFA} for a different proof). We will thus refer to representations of $\TE$ as `essential' or `singular' if the corresponding representations of $E$ are essential or singular, respectively.
 
If a representation of $E$ is essential, then the corresponding representation of $\TE$ will factor through the Cuntz-Pimsner algebra $\OO_E$ (see \cite{ms1}). However, not any representation of $\OO_E$ necessarily arises from an essential representation of $E$ (as a simple example, consider $E$ to be the correspondence over $\C$ of dimension $\aleph_0$; here, $\TE = \OO_E = \OO_{\infty}$, and if $s_1,s_2,...$ are the generators, and $\pi$ is the representation, then $\pi$ comes from an essential representation precisely when $\sum_{i=1}^{\infty}\pi(s_is_i^*) = 1$, with convergence taken in the strong operator topology). It seems, thus, natural to ask whether any $C^*$-correspondence as above admits an essential representation (we note that the non-degeneracy of the left action is surely needed).

More generally, we can ask the following question. Let $\pi_{ess}$ be the universal essential representation of $\TE$, i.e. a direct sum of unitary equivalence classes of all essential representations (where we bound the dimension of the target space $H$ to avoid set-theoretic difficulties).By the comments above, $\pi_{ess}$ factors through $\OO_E$. We preserve notation and think of $\pi_{ess}$ as a representation of $\OO_E$. The main result of this paper is the following.

\begin{Thm} 
\label{main theorem}
$\pi_{ess}$ is a faithful representation of $\OO_E$.
\end{Thm}

In particular, we see that $E$ always admits essential representations. 

This paper is structured as follows. In section \ref{decreasing functionals}, we adapt to our setting parts of Arveson's theory of decreasing weights (\cite{arveson 4, arveson's book}) and its discrete version due to Fowler (\cite{fowler}). The adaptation is straightforward, but we nontheless include full proofs for the reader's convenience. This will give us a sufficient criterion for a state on $\TE$ to give rise, via the GNS representation, to an essential representation of $E$.
In section \ref{construction}, we show how to construct sufficiently many such states to deduce Theorem \ref{main theorem}. We end by making some remarks on connections to endomorphisms.

I would like to thank Takeshi Katsura for a helpful discussion, and Michael Skeide for reading an earlier draft of this paper.

\section{Functionals on the algebra of operators of bounded support}
\label{decreasing functionals}

We fix some terminology about topologies on $\B(E)$. Let $T_{\lambda}$ be a net. We say that $T_{\lambda} \arrow 0$ in the \emph{strong topology} if $\|T_{\lambda}e\|_E \arrow 0$ for all $e \in E$, and in the \emph{strict topology} if $\|T_{\lambda}A\|_{\B(E)}, \|AT_{\lambda}\|_{\B(E)}\arrow 0$ for all $A \in \K(E)$ (i.e. the strict topology of $\B(E)$ coming from its being the multiplier algebra of $\K(E)$).
 Note that the strong topology is coarser than the strict topology.
We say that a linear functional on $\B(E)$ is strictly continuous if it is continuous with respect to the strict topology. Suppose $\varphi$ is a state on $\K(E)$, then we can (uniquely) extend $\varphi$ to a strictly continuous state on $\B(E)$ (using a GNS representation). 
Sesquilinear forms here will always be linear on the right, and conjugate-linear on the left. 

Let $P_n \in \B(\EE)$ be the projection onto the submodule $ \bigoplus_{k=0}^{n}E^{\otimes k}$. Denote $\J_n(\EE) = P_n\B(\EE)P_n$, and let $\Jc(\EE) = \bigcup_{n=0}^{\infty}\J_n(\EE)$ (no closure), the $*$-algebra of operators of bounded support.
 Let $Q_n = 1-P_{n-1}$, $R_n = P_n - P_{n-1}$ (i.e. $R_n$ is the projection onto $E^{\otimes n}$, and $Q_n$ is the projection onto $\bigoplus_{k=n}^{\infty}E^{\otimes k}$).
Let $\beta:B(\EE) \arrow B(Q_1\EE)$ be the representation given by $\beta(T) = T \otimes 1_E$, where we identify $Q_1\EE$ with $\EE \otimes_{\A} E$ (notice that to do so, we use the nondegeneracy of the left action). We think of $\beta$ as an endomorphism of $\B(\EE)$, where we identify $B(Q_1\EE)$ with $Q_1B(\EE) Q_1$. We now let $\delta: \B(\EE) \arrow \B(\EE)$ be the self-adjoint linear map given by $\delta(T) = T-\beta(T)$. Notice that $\beta(\Jc(\EE)) \subseteq \Jc(\EE)$ and therefore also $\delta(\Jc(\EE)) \subseteq \Jc(\EE)$.

Define an unbounded map $\lambda:\Jc(\EE) \arrow \B(\EE)$ by $\lambda(T) = \sum_{k=0}^{\infty}\beta^k(T)$, where the convergence is taken in the strong topology. Now, $\beta(\lambda(T)) = \lambda(T)-T$, so $\delta(\lambda(T)) = T$. Similarly, $\lambda(\delta(T)) = T$, so in particular, we have that $\lambda(\Jc(\EE)) \supseteq \Jc(\EE)$. $\lambda$ is clearly positive, being a sum of homomorphisms. 

Notice that we have 
\begin{equation}
\nonumber
\delta(xy) = x\delta(y)+\delta(x)y - \delta(x)\delta(y)
\end{equation}
Note that if $x \in \J_n(\EE)$ and $y\in\J_m(\EE)$ then $\lambda(x)y,x\lambda(y),xy \in J_{n+m}(\EE)$. 
It follows that $\lambda(x)\lambda(y) \in \lambda(\J_{n+m}(\EE))$. Therefore, $\lambda(\Jc(\EE))$ is a unital $*$-algebra.
Denote $\HTE =\overline{\lambda(\Jc(\EE))}$. Notice that $\beta$ preserves $\HTE$. We shall denote the restriction of $\beta$ to this $C^*$-algebra by $\beta$ as well.

Let $\xi \in E^{\otimes n}$, $\eta \in E^{\otimes m}$, then $\beta(T_{\xi} T_{\eta}^*) \zeta = T_{\xi}T_{\eta}^*\otimes 1_E \zeta = 0$ if $\zeta \in E^{\otimes k}$, $k \leq m$, and otherwise, we have $\beta(T_{\xi} T_{\eta}^*) \zeta = T_{\xi}T_{\eta}^* \zeta$. Consequently, we have that $\delta (T_{\xi} T_{\eta}^*) = \xi \otimes \eta^* \in \K(\EE)$. Therefore, we have that $\delta(\TE) \subseteq \K(\EE)$, and $\lambda (\Jc(\EE) \cap \K(\EE))$ is dense in $\TE$. 

\begin{Rmk} The fact that $\delta(\TE) \subseteq \K(E)$ appears in \cite{pimsner}, where it is used to obtain an element of $KK(\TE,\A)$ (shown to be invertible) from the pair $(id,\beta)$ (in our notation).
\end{Rmk}

\begin{Def} A linear functional $\omega$ on $\Jc(\EE)$ is said to be \emph{decreasing} if $\omega \circ \delta$ is a positive functional on $\Jc(\EE)$. 
\end{Def}

We recall a simple technical lemma characterizing when positive linear functionals on a normed algebra are bounded (see \cite{arveson's book}, lemma 4.9.4). We include a proof for the reader's convenience, but just state it for the unital case, which is the one we need.

\begin{Lemma} 
\label{boundedness condition}
Let $\A$ be a unital normed complex algebra with isometric involution. Let $\varphi$ be a linear functional on $\A$, which is positive, in the sense that $\varphi(x^*x) \geq 0$ for all $x \in \A$. If $\limsup_{n \arrow \infty}|\varphi(x^n)|^{1/n} \leq \|x\|$ for all $x$, then $\varphi$ is bounded, with norm $\varphi(1)$.
\end{Lemma}
\begin{proof}
By the Cauchy-Schwarz inequality, $|\varphi(x)| = |\varphi(1 \cdot x)| \leq (\varphi(1)\varphi(x^*x))^{1/2}$. Applying this again to the right hand side $n$ times yields 
$$
|\varphi(x)| \leq \varphi(1)^{1/2 + 1/4 + ... + 1/2^n}\varphi\left ( (x^*x)^{2^{n-1}} \right )^{1/2^n}
$$ 
and by taking $\limsup$, we get that $|\varphi(x)| \leq \varphi(1)\|x^*x\|^{1/2} \leq \varphi(1)\|x\|$, as required.
\end{proof}

The following proposition is a straightforward generalization of \cite{fowler}, Proposition 1.5. 

\begin{Prop} 
\label{positive iff positive and decreasing}
$\omega \circ \delta$ is positive on $\lambda(\Jc(\EE))$ if and only if $\omega$ is positive and decreasing (i.e. if both $\omega$ and $\omega \circ \delta$ are positive on $\Jc(\EE)$) in which case $\|\omega \circ \delta\| = \omega(P_0)$. 
\end{Prop}
\begin{proof}
If $\omega \circ \delta$ is positive on $\lambda(\Jc(\EE))$, then it is also positive on $\Jc(\EE) \subseteq \lambda(\Jc(\EE))$, and $\omega = \omega \circ \delta \circ \lambda$ is positive, as it is a composition of positive functionals.

Now, suppose $\omega$ is positive and decreasing. Since $\omega$ is positive, we can use it to view $\Jc(\EE)$ as an semi-inner product space, and complete it to a Hilbert space $H$. Denote by $\hat{x}$ the image of $x \in \Jc(\EE)$ in $H$ (i.e. the coset of $x$ module the null-space of $\omega$ in the completion). $\omega$ is decreasing, so $\omega(x^*x) - \omega(\beta(x^*x)) \geq 0$ for all $x \in \Jc(\EE)$. Therefore, the map $\hat{x} \mapsto \widehat{\beta(x)}$ is well defined, and extends (uniquely) to a linear contraction $T:H \arrow H$.  

Set 
$$
T_k = \left \{ \begin{matrix}T^k &|& k \geq 0 \\ T^{*-k}
 &|& k<0 \end{matrix} \right .
$$
Seen as operator valued function $\Z \arrow \B(H)$, the map $k \mapsto T_k$ is positive definite (\cite{sz-nagy-foias}).
One checks that for $x\in \Jc(\EE)$, we have 
$$
\delta(\lambda(x)^*\lambda(x)) = x^*\lambda(x)+\lambda(x^*)x-x^*x = \sum_{n=0}^{\infty}x^*\beta^n(x) + \sum_{n=1}^{\infty}\beta^n(x^*)x
$$ 
(a finite sum). Notice that 
$$
\lb \hat{x},T_k\hat{x}\rb_H =  \left \{ \begin{matrix} \omega(x^*\beta^k(x)) &|&k \geq 0 \\ \omega(\beta^{-k}(x^*)x) &|& k < 0\end{matrix} \right .
$$ 
So, we have $\omega \circ \delta (\lambda(x)^*\lambda(x)) = \sum_{k=-\infty}^{\infty}\lb \hat{x},T_k\hat{x} \rb \geq 0$ (being the value of the Fourier transform at $z=1$ of the finitely supported positive definite function on $\Z$ given by $k \mapsto \lb \hat{x},T_k\hat{x} \rb$).

It remains to prove that $\omega \circ \delta$ is bounded. For that, we use lemma \ref{boundedness condition}. Notice that $\omega \circ \delta (1) = \omega (P_0)$. So, it remains to show that for any $k$, and any $x \in \J_k(\EE)$, we have $\limsup_{n \arrow \infty}|\omega\circ\delta(\lambda(x)^n)|^{1/n} \leq \|\lambda(x)\|$. Notice that $\lambda(x)^n \in \lambda(\J_{kn}(\EE))$, and that $\|\delta\| = 2$. So, $\|\omega \circ \delta (\lambda(x)^n)\| \leq \|\omega|_{\J_{kn}(\EE)}\|\|\delta((\lambda(x)^n)\| \leq 2 \omega(P_{k+n-1})\|\lambda(x)\|^n$. Now, $P_{kn} = \sum_{j=0}^{kn}\beta^j(P_0)$. Since $\omega$ is decreasing, we have then that $\omega(P_{kn}) \leq (kn+1)\omega(P_0)$. So, $|\omega\circ\delta(\lambda(x)^n)|^{1/n} \leq \left ( 2(kn+1) \omega (P_0) \right )^{1/n} \|\lambda(x)\|$, which gives us the required inequality.
\end{proof}

Slightly abusing notation, we denote the (necessarily unique) extention of $\omega \circ \delta$ to $\HTE$ by $\omega \circ \delta$ as well. Notice that if $\varphi$ is a positive linear functional on $\HTE$, then $\omega = \varphi \circ \lambda$ is positive and decreasing, so any positive linear functional on $\HTE$ is obtained in this way.

A linear functional $\omega$ on $\Jc(\EE)$ is said to be \emph{locally strict} if for all $n$, the restiction of $\omega$ to $\J_n(\EE)$ is strictly continuous.

Suppose $\omega_1,\omega_2$ are two locally strict decreasing functionals on $\Jc(\EE)$, and the restrictions of $\omega_i \circ \delta$ $i= 1,2$ agree on $\TE = \overline{\lambda (\Jc(\EE) \cap \K(\EE))}$, then $\omega_1$ and $\omega_2$ agree on $\Jc(\EE) \cap \K(\EE)$, and since they are locally strict, $\omega_1,\omega_2$ agree everywhere. 

Suppose now that $\varphi$ is a positive linear functional on $\TE$. Define $\omega$ on $\Jc(\EE) \cap \K(\EE)$ by $\omega = \varphi \circ \lambda$. $\omega$ is positive (and bounded) on $\J_n(\EE) \cap \K(\EE)$, and therefore extends, uniquely, to a strictly continuous functional on $\J_n(\EE)$. We view $\omega$ as a locally strict functional on $\J_c(\EE)$. 

\begin{Prop} $\omega$, obtained as described in the previous paragraph, is positive and decreasing.
\end{Prop}
The proof is an easy generalization of that of \cite{fowler}, Theorem 1.9.
\begin{proof}
It suffices to show that $\omega \circ \delta$ is positive on $\lambda(\J_c(\EE))$. Notice that we already know that $\omega \circ \delta$ is positive on $\lambda(\J_c(\EE) \cap \K(\EE))$. For any $x \in \J_c(\EE)$, define a function $f_x:\Z \arrow \C$ by 
$$
f_x(k) = \left \{ \begin{matrix} \omega(x^*\beta^k(x)) &|&k \geq 0 \\ \omega(\beta^{-k}(x^*)x) &|& k < 0\end{matrix} \right .
$$
As in the proof of proposition \ref{positive iff positive and decreasing}, we see that $\omega \circ \delta (\lambda(x)^*\lambda(x)) = \sum_{k=-\infty}^{\infty}f_x(k)$ (only finitely summands are non-zero), so it will suffice if we could show that $f_x$ is in fact positive definite. Since $\omega$ is locally strict, $f_x$ is a pointwise limit of functions of the form $f_y$, where $y \in \J_c(\EE) \cap \K(\EE)$. It therefore suffices to show that $f_x$ is positive definite when $x \in \J_c(\EE) \cap \K(\EE)$, and therefore, it suffices to show that the Fourier transform $\hat{f}_x(z) \geq 0$ for all $z \in \T$. 

Suppose $x\in \J_N(\EE)\cap \K(\EE)$. For $z \in \T$, let $u_z = \sum_{n=0}^{N} z^nR_n$. Notice that $u_zx \in \J_N(\EE) \cap \K(\EE)$, and a simple computation shows that $f_{u_zx}(k) = \bar{z}^kf_x(k)$. So, $\hat{f}_x(z) = \sum_{k=-\infty}^{\infty} f_x(k)\bar{z}^k =  \sum_{k=-\infty}^{\infty} f_{u_zx}(k) = \omega\circ\delta(\lambda(u_zx)^*\lambda(u_zx)) \geq 0$, as required.
\end{proof}

Thus, for any positive linear functional $\varphi$ on $\TE$ we have a positive linear functional $\hat{\varphi}$ on $\HTE$, given by taking $\omega = \varphi \circ \lambda$, extending it to $\J_c(\EE)$, and composing with $\delta$. Note that any approximate unit of $\A$ is also an approximate unit of $\TE$ (when we view $\A$ as a subalgebra of $\TE$), since we assumed that the left action of $\A$ on $E$ is non-degenerate. Therefore, any state $\varphi$ of $\TE$ restricts to a state of $\A$. In the above notation, we have that $\omega(P_0) = \|\varphi|_{\A}\|$, so, if $\varphi$ is a state, then so is $\hat{\varphi}$. 

\begin{Rmk} The map $S(\TE) \arrow S(\HTE)$ given by $\varphi \mapsto \hat{\varphi}$ is clearly affine, although the reader is cautioned that generally it is not continuous.
\end{Rmk}

\begin{Def} We call a state of $\TE$ \emph{essential} if its associated GNS representation, restriced to the canonical copies of $\A$,$E$ in $\TE$, is essential.
\end{Def}

Let $\varphi$ be a state of $\TE$, and let $\omega$ be its associated locally strict functional.

\begin{Claim} If $\varphi = \hat{\varphi} \circ \beta$ then $\varphi$ is essential.
\label{fc characterization}
\end{Claim}
\begin{proof}
Let $(\pi,H,\Omega)$ be the GNS representation of $\TE$ with respect to $\varphi$.  Let $Q'$ be the projection onto $\bigcap_{n=0}^{\infty}\overline{\pi(E^{\otimes n})H}$, then $Q'H$ is invariant, the restriction of $\pi$ to $Q'H$ is essential, and $Q'H$ is the largest subspace of $H$ which satisfies this property. Thus, $\varphi$ is essential if and only if $Q'\Omega = \Omega$, i.e. if $\|Q'\Omega\| = \lb \Omega, Q'\Omega \rb = 1$. Let $Q'_n = \overline{\pi(E^{\otimes n})H}$ (so $Q'_n \arrow Q'$ in the strong operator topology, and so, $\lb \Omega, Q'_n\Omega \rb \arrow \lb \Omega, Q'\Omega \rb$)). If $a_k^{(n)}$ is an approximate identity for $\K(E^{\otimes n})$ then, $\pi \left ( \lambda(a_k^{(n)}) \right ) \arrow Q'_n$ in the strong operator topology. Therefore,  $\lb \Omega, Q'_n\Omega \rb = \lim_k \lb \Omega, \lambda(a_k^{(n)})\Omega \rb = \lim_k \varphi (\lambda(a_k^{(n)})) = \lim_k \omega (a_k^{(n)}) = \omega(R_n) = \hat{\varphi}(Q_n) = \hat{\varphi}\circ \beta^n(1)$. 
Thus, if $\varphi = \hat{\varphi} \circ \beta$, then, if $(a_k)$ is an approximate unit for $\A$, then $\hat{\varphi}\circ\beta^n(1) = \lim_k \hat{\varphi}\circ\beta^n(a_k) = \lim_k\varphi(a_k) = 1$  and so $\varphi$ is essential.
\end{proof}

\begin{Claim} $\varphi = \hat{\varphi} \circ \beta$ if and only if $\omega = \omega \circ \beta$.
\label{fc characterization 2}
\end{Claim}
\begin{proof} Note that $\lambda \circ \beta = \beta \circ \lambda$, and  $\delta \circ \beta = \beta \circ \delta$. So if $\varphi = \hat{\varphi} \circ \beta$ then $\omega = \hat{\varphi} \circ \lambda$, and $\omega  \circ \beta= \hat{\varphi} \circ \lambda \circ \beta = \hat{\varphi} \circ \beta \circ \lambda = \varphi \circ \lambda = \omega$. Conversely, if $\omega = \omega \circ \beta$ then $\varphi = \omega \circ \delta =  \omega \circ \beta \circ \delta = \hat{\varphi}\circ \beta$.
\end{proof}
  
\section{Construction of essential states}
\label{construction}

By Claims \ref{fc characterization} and \ref{fc characterization 2}, we can construct essential states on $\TE$ as follows. Suppose we have a sequence of strictly continuous states $\nu_n$ on $\B(E^{\otimes n})$, such that $\nu_{n+1} \circ \beta = \nu_n$ for all $n$. Viewing each such $\nu_n$ as a positive linear functional on $\J_c(\EE)$ (with support $R_n$), we can form $\omega = \sum_n\nu_n$. $\omega$, then, is clearly a locally strict decreasing functional, and $\omega \circ \beta = \omega$, so the associated state $\varphi$ on $\TE$ is essential.

Let $\nu_0$ be an arbitrary state on $\A$ (extended as strictly continuous state to $\M(\A) \cong \B(E^{\otimes 0})$ if $\A$ is non-unital). We can push forward $\nu_0$ to $\beta(\B(E^{\otimes 0}))$, thought of as a $C^*$-subalgebra of $\B(E)$, and extend it to a state $\nu_1$ on $\B(E)$. Now push $\nu_1$ to the subalgebra $\beta(\B(E))$ of $\B(E^{\otimes 2})$, etc. The problem is that while we can always extend the state $\nu_n$, thought of as a state of $\beta(\B(E^{\otimes n}))$ to all of $\B(E^{\otimes n+1})$, it is not necessarily the case that the extension can be chosen to be strictly continuous, as the following example shows.

\begin{Exl}
Let $\A = C([0,1])$, $E = C([0,1)) \oplus C((0,1])$ (with the standard inner product), and take the left action to be 
$$
f \cdot (g_1,g_2) (x) = (f(x/2)g_1(x),f((x+1)/2)g_2(x))
$$ 
Any state $\nu_0$ on $\A$ is given by a measure on $[0,1]$, and it is easy to check that we can choose a strictly continuous extension to $\B(E)$ if and only if $1/2$ is not an atom of this measure.
\end{Exl}

We will show, though, that there is a dense set of states of $\A$ which do admit a strictly continuous sequence of extensions as above. To that end, we start with some technical lemmas.

\begin{Lemma} 
\label{Baire category}
Let 
$$
\begin{CD} X_0 @<f_{0,1}<<X_1 @<f_{1,2}<< X_2 @<f_{2,3}<< \cdots \end{CD}
$$  
be an inverse system of surjective continuous maps, where $X_n$, $n=0,1,...$ are compact Hausdorff spaces. Denote by $f_{n,m}:X_m \arrow X_n$, $n<m$ be the map obtained by composition (and $f_{n,n}$ will be the identity). Let $X_{\infty}$ be the inverse limit, and let $f_{n,\infty}:X_{\infty} \arrow X_n$ be the canonical maps. 

Suppose we have dense open sets $U_n \subseteq X_n$, $n = 0,1,2,...$ such that $f_{n,n+1}(U_{n+1}) \subseteq U_{n}$ for all $n$, then $f^{-1}_{n,\infty}(U_n)$ are dense open sets in $X_{\infty}$. 

In particular, $\bigcap_{n=0}^{\infty}f^{-1}_{n,\infty}(U_n)$ is dense in $X_{\infty}$. Thus, there exists a dense collection of $x_0\in X_0$ for which there is a sequence $x_n \in U_n$ $n=1,2...$ with $f_{n,m}(x_m) = x_n$ for all $n<m$.
\end{Lemma}
\begin{proof}
The only thing which might not seem immediate is that $f^{-1}_{n,\infty}(U_n)$ are dense. So, we need to show that any non-empty open set in $X_{\infty}$, has non-empty intersection with $f^{-1}_{n,\infty}(U_n)$. It suffices to do so for basic open sets: let $V_k \subseteq X_k$, $k=0,...,m$ (for an arbitrary $m$, which we may assume is $\geq n$) be open sets such that $\bigcap_{k=0}^m f^{-1}_{k,m}(V_k) \neq \emptyset$. It suffices to show that $f^{-1}_{n,\infty}(U_n) \cap \bigcap_{k=0}^m f^{-1}_{k,\infty}(V_k) \neq \emptyset$. So, 
$$f^{-1}_{n,\infty}(U_n) \cap \bigcap_{k=0}^m f^{-1}_{k,\infty}(V_k) \supseteq 
f^{-1}_{m,\infty}(U_m) \cap \bigcap_{k=0}^m f^{-1}_{m,\infty} \left ( f^{-1}_{k,m} (V_k) \right )  = 
f^{-1}_{m,\infty} \left ( U_m \cap \bigcap_{k=0}^m f^{-1}_{k,m} (V_k) \right )
$$ 
which is nonempty since $U_m$ is dense in $X_m$ and $\bigcap_{k=0}^m f^{-1}_{k,m} (V_k)$ is open and non-empty.
\end{proof}

\begin{Lemma} 
\label{density}
Let $\B$ be a $C^*$-algebra. The strictly continuous states on $\M(B)$ are a co-meager (dense $G_{\delta}$) set in the state space of $\M(B)$. 
\end{Lemma}
\begin{proof} Consider the natural embedding of $\B$ in $\B^{**}$ (the double dual). The states of $\B$ (which we identify with the strictly continuous states on $\M(B)$, as well as the normal states on $\B^{**}$) are dense in the state space of $\B^{**}$, and therefore also in the state space of $\M(\B)$. Let $\{e_{\lambda}\}_{\lambda\in\Lambda}$ be an increasing approximate identity for $\B$. A state $\varphi$ on $\M(\B)$ is strictly continuous if and only if $\sup_{\lambda\in\Lambda}\varphi(e_{\lambda}) = 1$. So, $U_m = \{\varphi \in S(\M(\B)) \; | \; \exists \lambda\in\Lambda \textrm{ s.t. }\varphi(e_{\lambda}) > 1-1/m\}$ is open (and dense, since it contains the strictly continuous states), so the strictly continuous states are the intersection of open dense sets $\bigcap_m U_m$, as required. 
\end{proof}

\begin{Lemma} 
\label{norm of restriction}
Suppose $\A$ is a subalgebra of $\M(\B)$ which is non-degenerate, in the sense that $\A\B = \B$. If $\varphi$ is a state on $\M(\B)$, then $\|\varphi|_{\A}\| \geq \|\varphi|_{\B}\|$. In particular, if $\varphi$ is a 
strictly continuous state of $\M(\B)$, then the restriction of $\varphi$ to $\A$ is a state of $\A$. 
\end{Lemma}
\begin{proof} 
Let $a_{\lambda}$ be an approximate unit for $\A$, then non-degeneracy implies that for any $b \in \B$, $\|a_{\lambda}b-b\| \arrow 0$. Let $r =  \|\varphi|_{\B}\|$. Given $\epsilon>0$, choose $b\in \B$ positive of norm 1 such that $\varphi(b)>r-\epsilon/2$, and $\lambda$ such that $ \|a_{\lambda}b-b\| < \epsilon/2$, so $r-\epsilon < \varphi(b)-\epsilon/2 < |\varphi(a_{\lambda}b)| \leq \varphi(a_{\lambda}^*a_{\lambda})\varphi(b^*b) \leq \varphi(a_{\lambda}^*a_{\lambda})$, so the restriction of $\varphi$ to $\A$ has norm at least $r$, as required.
\end{proof}

\begin{Lemma} Let $\A_0,\A_1,...$ be $C^*$-algebras, and suppose we are given non-degenerate injective homomorphisms $\A_n \arrow \M(\A_{n+1})$, $n=0,1,...$, that is, we have a diagram as follows (where the map $\A_n \arrow \M(\A_n)$ is the inclusion).
$$
\xymatrix{\M(\A_0) \ar[r]^{\varphi_{1,0}} & \M(\A_1) \ar[r]^{\varphi_{2,1}} & \M(\A_2) \ar[r]^{\varphi_{3,2}} & \cdots \\
             \A_0 \ar[u] \ar[ur] & \A_1 \ar[u] \ar[ur] & \A_2 \ar[u] \ar[ur] & \cdots
}$$
We use the same notation for a state on $\A_n$ and its strictly continuous extension to $\M(\A_n)$. It follows that there is a dense set of states $S_0$ of $\A_0$, such that each $\nu_0 \in S_0$ admits a sequence of strictly continuous states $\nu_n \in S(\M(\A_n))$ with $\nu_{n+1} \circ \varphi_{n+1,n} = \nu_n$.  
\end{Lemma}
\begin{proof} Denote by $\varphi_{m,n}:\M(\A_n) \arrow \M(\A_m)$, $n<m$ the map obtained by composition, and by $f_{n,m} = \varphi_{m,n}^* : S(\M(\A_m)) \arrow S(\M(\A_n))$ the induced map on the state space.
As in the proof of Lemma \ref{density}, set 
$$
U_m = \{\varphi \in S(\M(\A_m)) \; | \; \exists \lambda\in\Lambda \textrm{ s.t. }\varphi(e_{\lambda}) > 1-1/(m+1)\}
$$ 
where $\{e_{\lambda}\}_{\lambda \in \Lambda}$ is some approximate unit for $\A_m$. $U_m$ is open and dense in $S(\M(\A_m))$, and we have $f_{n,n+1}(\nu) \in U_n$ (by Lemma \ref{norm of restriction}) for all $\nu \in U_{n+1}$. It follows now from Lemma \ref{Baire category} that there are states $\nu_n \in U_n$, $n=0,1,2,...$ such that $f_{n,m}(\nu_m) = \nu_n$ for all $n<m$. Since $\nu_n \in \bigcap_{m \geq n} f_{n,m}(U_m)$, it follows that $\nu_n$ must be strictly continuous, as required. 
\end{proof}

\begin{Cor} There is a dense set of states $\nu_0$ of $\A$ which admit a sequence of strictly continuous exentsions $\nu_n \in S(\B(E^{\otimes n}))$, as described in the beginning of the section.
\end{Cor}
\begin{proof} Apply the previous lemma to 
$$
\xymatrix{\M(\A) \ar[r] & \B(E) \ar[r] & \B(E^{\otimes 2}) \ar[r] & \cdots \\
             \A \ar[u] \ar[ur] & \K(E) \ar[u] \ar[ur] & \K(E^{\otimes 2}) \ar[u] \ar[ur] & \cdots
}$$
\end{proof}

Let $\psi$ be a representation of $E$ on $H$. For $z\in \T$, we can define a new representation $\psi_z$ of $E$ on $H$ given by $\psi_z(a)  = \psi(a)$ ($a \in \A$), and $\psi_z(e) = z\psi(e)$ ($e \in E$). Clearly, if $\psi$ is essential then so is $\psi_z$. Recall the following theorem, from \cite{FMR} (section 4). Let $\gamma_z$ be the automorphism of $\OO_E$ given by mapping $e \mapsto ze$, and fixing $\A$ (where we view $\A$ and $E$ as contained in $\OO_E$).

\begin{Thm}[Gauge invariant uniqueness theorem] 
\label{GIUT}
Let $\pi$ be a representation of $\OO_E$. If there is a point-norm continuous action $\alpha_z : \T \arrow Aut(\pi(\OO_E))$ such that $\alpha_z \circ \pi = \pi \circ \gamma_z$ for all $z\in \T$, and if $\pi|_{\A}$ is faithful, then $\pi$ is faithful.
\end{Thm}

We can now prove the main theorem.

\begin{proof}[Proof of Theorem \ref{main theorem}] 
If $\varphi$ is a state of $\TE$, and $\omega$ is its corresponding decreasing functional, and $a \in \A$, then $\varphi(a) = \omega(a)$, where $\A$ is thought of as a subalgebra of $\TE$ and as $\J_0(\EE)$, respectively. Let $S$ be the set of states $\nu_0$ of $\A$ which admit strictly continuous extensions $\nu_1,\nu_2,...$ as above. Since $S$ is dense, it follows that for any $a \in \A_+$ and $\epsilon>0$ there is an essential representation $\psi$ of $E$ such that $\|\psi(a)\| > \|a\|-\epsilon$ (we pick $\nu_0 \in S$ such that $\nu_0(a)  > \|a\|-\epsilon$, let $\varphi$ be the corresponding state of $\TE$, let $\pi_{\varphi}$ be the GNS representation, and take $\psi$ to be the representation of $E$ associated to $\pi_{\varphi}$). Thus, if we let $\pi_{ess}$ be the universal essential representation as above, we see that $\pi_{ess}|_{\A}$ is faithful. By the remarks preceeding Theorem \ref{GIUT}, $\pi_{ess}$ satisfies the hypothesis of the theorem, and therefore, $\pi_{ess}$ is faithful, as required. 
\end{proof}

\section{Concluding remarks}

Any representation $\psi$ of $E$ on $H$ gives rise to an endomorphism $\alpha$ of $\psi(\A)'$, by $\alpha(x)\psi(e)\xi = \psi(e)x\xi$ (for $x \in \psi(\A)'$, $e \in E$, $\xi \in H$, and defined to be 0 on $\psi(E)H^{\perp}$; the reader is referred to \cite{ms1} for further details, including why this is well defined). $\alpha$ is unital if and only if $\psi$ is essential. 

Now, if $\M$ is a von-Neumann subalgebra of $\B(H)$ and $\alpha$ is a normal endomorphism of $\M'$, then $E_{\alpha} = \{T \in \B(H) \; | \; \alpha(x)T = Tx\}$ is naturally a $W^*$-correspondence over $\M$, where by a $W^*$-correspondence we mean a $C^*$-correspondence $E$ over a von-Neumann algebra $\M$, where $E$ is self-dual (see \cite{paschke}), and the left action $\M \arrow \B(E)$ is normal (recall that in this case, $\B(E)$ is naturally a von-Neumann algebra). $\alpha$ is unital if and only if this concrete representatation of $E_{\alpha}$ is essential.
It is thus interesting to know whether any $W^*$-correspondence arises as the intertwining space of a unital normal endomorphism. Recall that if $E$ is a self-dual Hilbert module over $\M$, then $E$ has a predual, which gives it a weak$^*$-topology (which we'll call the ultraweak topology). Thus, the concrete question is the following.
\begin{Question} Let $E$ be a $W^*$-correspondence. Does $E$ admit an essential representation $\psi$, such that $\psi_{\M}$ and $\psi_{E}$ are continuous when $E$, $\M$ and $\B(H)$ are endowed with their ultraweak topologies?
\end{Question}
Our construction guarantees the existence of essential representations, however it is unclear if they are ultraweakly continuous. It is worth noting that in this setting, if $\nu_0$ is a normal state of $\M$, then we can obtain a sequence of strictly continuous states $\nu_n \in \B(E^{\otimes n})$ as discussed in section \ref{construction}. Thus, the corresponding state $\varphi$ of the $C^*$-algebra $\TE$ will satisfty that its restriction to $\M \subset \TE$ be ultraweakly continuous. However, this in its own does not appear to guarantee that the restriction of the associated GNS representation to $\M$ be ultraweakly continuous. 
Of course, if $\A$ is finite dimensional, then ultraweak continuity is immediate. 

Much of the motivation for this present work comes from Arveson's work on $E_0$-semigroups, i.e. continuous one-parameter semigroups of unital endomorphisms of $\B(H)$ (see \cite{arveson's book}). If $\alpha$ is an $E_0$-semigroup of $\B(H)$, the intertwiners of each $\alpha_x$, $x>0$, form a Hilbert space in $\B(H)$, and together, they form a so-called product system of Hilbert spaces. Much of this generalizes to $E_0$-semigroups of more general von-Neumann algebras (and certain $C^*$-algebras), by replacing Hilbert spaces by correspondences. In \cite{arveson 4}, Arveson showed that every product system of Hilbert spaces is associated to an $E_0$-semigroup in this fashion, by constructing essential representations of product systems. 
There is hope that at least for the finite dimensional case, it might be possible to extend Arveson's result, by constructing a continuous analogue of the argument herein for certain continuous analogues of Pimsner's Toeplitz algebras (considered in \cite{me Crelle, me JFA}), however this goes beyond the scope of this paper.
We refer the reader to \cite{ms2,skeide} and references therein for more information on product systems of Hilbert modules in the context of $E_0$-semigroups.

\end{document}